\documentclass[11pt]{article}
\usepackage{amsmath,amssymb}

\newtheorem{theo}{Theorem}[section]
\newtheorem{prop}[theo]{Proposition}
\newtheorem{lemma}[theo]{Lemma}

\newcommand\npf{\mbox{ }\hfill\sqr\vskip6pt}
\def\sqr{$\vcenter{\hrule height.2mm
\hbox{\vrule width.2mm height2mm\kern2mm
\vrule width.2mm}\hrule height.2mm}$}

\begin{document}

\title{The order of monochromatic subgraphs with a given minimum degree}

\author{
Yair Caro
\thanks{Department of Mathematics, University of Haifa at Oranim, Tivon 36006, Israel. e--mail: yairc@macam98.ac.il}
\and
Raphael Yuster
\thanks{Department of Mathematics, University of Haifa at Oranim, Tivon 36006, Israel. e--mail: raphy@macam98.ac.il}
}

\date{}

\maketitle

\begin{abstract}
Let $G$ be a graph. For a given positive integer $d$, let $f_G(d)$ denote the largest
integer $t$ such that in every coloring of the edges of $G$ with two colors there is a monochromatic
subgraph with minimum degree at least $d$ and order at least $t$.
For $n > k > d$ let $f(n,k,d)$ denote the minimum of $f_G(d)$ where $G$ ranges over all
graphs with $n$ vertices and minimum degree at least $k$.
In this paper we establish $f(n,k,d)$ whenever $k$ or $n-k$ are fixed, and $n$ is sufficiently large.
We also consider the case where more than two colors are allowed.
\end{abstract}

\setcounter{page}{1}

\section{Introduction}
All graphs considered in this paper are finite, simple and undirected. For standard terminology
used in this paper see \cite{We}.
It is well known that in any coloring of the edges of a complete graph with two colors there is a monochromatic
connected spanning subgraph. This folkloristic Ramsey-type fact has been generalized in many ways,
where one shows that some given properties of a graph $G$ suffice in order to guarantee a large monochromatic subgraph
of $G$ with related given properties in any two (or more than two) edge-coloring of $G$.
See, e.g., \cite{BiDiVo,Ma,CaRo,CaYu} for these types of results.
In this paper we consider the property of having a certain minimum degree.
For given positive integers $d$ and $r$, let $f_G(d,r)$ denote the largest
integer $t$ such that in every coloring of the edges of the graph $G$ with $r$ colors there is a monochromatic
subgraph with minimum degree at least $d$ and order at least $t$.
For $n > k > d$ let $f(n,k,d,r)$ denote the minimum of $f_G(d)$ where $G$ ranges over all
graphs with $n$ vertices and minimum degree at least $k$.
The main results of our paper establish $f(n,k,d,2)$ whenever $k$ or $n-k$ are fixed, and $n$ is sufficiently large.
In particular, we prove the following results.
\begin{theo}
\label{t1}
(i) For all $d \geq 1$ and $k \geq 4d-3$,
\begin{equation}
\label{e1}
f(n,k,d,2) \geq \frac{k-4d+4}{2(k-3d+3)}n+ \frac{3d(d-1)}{4(k-3d+3)}.
\end{equation}

\noindent
(ii) For all $d \geq 1$ and $k \leq 4d-4$, if $n$ is sufficiently large then $f(n,k,d,2) \leq d^2-d+1$.
In particular, $f(n,k,d,2)$ is independent of $n$.
\end{theo}
\begin{theo}
\label{t2}
For all $d \geq 1$, $r \geq 2$ and $k > 2r(d-1)$, there exists an absolute constant $C$ such that
$$
f(n,k,d,r) \leq n \frac{k-2r(d-1)}{r(k-(r+1)(d-1))}+C.
$$
In particular, $f(n,k,d,2) \leq \frac{k-4d+4}{2(k-3d+3)}n+C$.
\end{theo}
Notice that Theorem \ref{t1} and Theorem \ref{t2} show that for fixed $k$, $f(n,k,d,2)$ is determined up to
a constant additive term.
The following theorem determines $f(n,k,d,2)$ whenever $k$ is very close to $n$.
\begin{theo}
\label{t3}
Let $d$ and $k$ be positive integers. For $n$ sufficiently large, $f(n,n-k,d,2) = n -2d-k+3$.
\end{theo}
The next section presents our main results. The final section contains some concluding remarks.
Throughout the rest of this paper, we use the term {\em $k$-subgraph} to denote a subgraph with minimum
degree at least $k$. 

\section{Results}
We need the following lemmas.
The first one is well-known (see, e.g., \cite{Bo} page xvii).
\begin{lemma}
\label{l21}
For every $m \geq k$, every graph with $m$ vertices and more than $(k-1)m-{k \choose 2}$
edges contains a $k$-subgraph. Furthermore, there are graphs with $m$ vertices and
$(k-1)m-{k \choose 2}$ edges that have no $k$-subgraph. \npf
\end{lemma}
\begin{lemma}
\label{l22}
Let $X$ be the set of at least $k$ vertices of a graph $G$ that are not on any $k$-subgraph.
Then, the sum of the degrees of the vertices of $X$ is at most $2(k-1)|X|-{k \choose 2}$.
\end{lemma}
{\bf Proof}\,
Assume the lemma is false.
Put $x=|X|$ and let $S \subset V(G) \setminus X$ denote the set of vertices of the graph $G$
that have a neighbor in $X$. Put $s=|S|$. Notice that there are at most $(k-1)x-{k \choose 2}$ edges
with both endpoints in $X$, and hence, if $z$ denote the number of edges between $X$ and $S$ then $z > {k \choose 2}$.
We distinguish between two cases. Assume first that $s \geq k$.
Replace the edges of $G$ with both endpoints in $S$ with a set $M$ of
$(k-1)s-{k \choose 2}$ edges that induce no $k$-subgraph (such an $M$ exists by Lemma \ref{l21}).
After this replacement, the sum of the degrees of the subgraph on $X \cup S$
is greater than
$$
2(k-1)x-2{k \choose 2}+2z+2(k-1)s-2{k \choose 2} \geq 2(k-1)(x+s)-k(k-1).
$$
Hence, this subgraph which has $x+s$ vertices, has more than $(k-1)(x+s)-{k \choose 2}$
edges and therefore contain a $k$-subgraph, $P$. Clearly, $P$ contains at least one vertex
of $X$. Now, delete $M$ and restore the original edges with both endpoints in $S$.
Also, add to $P$ all other vertices of $V(G) \setminus (X \cup S)$ and all their
incident edges. The obtained graph is a $k$-subgraph of $G$ that contains a vertex of $X$,
a contradiction. Now assume $s < k$ (clearly $s \geq 1$). We can repeat the same argument where instead of $M$
we use a complete graph on $S$, and similar computations hold.
\npf

\noindent
{\bf Proof of Theorem \ref{t1}, part (i).}\,
The theorem is trivial for $d=1$ so we assume $d \geq 2$.
Let $G=(V,E)$ have $n$ vertices and minimum degree at least $k$,
and consider some fixed red-blue coloring of $G$.
Let $B$ ($R$) denote the set of vertices of $G$ that are not on any blue (red) $d$-subgraph
but are on some red (blue) $d$-subgraph.
Let $C$ denote the set of vertices that are not on any red nor blue $d$-subgraph.
Put $|R|=r$, $|B|=b$, $|C|=c$. Clearly, there is a monochromatic subgraph of order at
least $(n-|C|)/2$. Hence, if $|C| < d$ the theorem trivially holds since the r.h.s. of (\ref{e1}) is always
at most $(n-d+1)/2$. We may therefore assume $|C| \geq d$.
For each $v \in B \cup C$ ($v \in R \cup C$) let $b(v)$ ($r(v)$) denote
the number of blue (red) edges incident with $v$ and that are not on any blue (red) $d$-subgraph.
By Lemma \ref{l22},
$$
\sum_{v \in B \cup C} b(v) \leq 2(d-1)(b+c)-{d \choose 2}, \qquad \sum_{v \in R \cup C} r(v) \leq 2(d-1)(r+c)-{d \choose 2}.
$$
Notice that, trivially, for each $v \in C$, $b(v)+r(v)=deg(v) \geq k$.
Put
$$
b_c = \sum_{v \in C}b(v), \qquad r_c = \sum_{v \in C}r(v).
$$
Thus, $b_c+r_c \geq kc$.
By Lemma \ref{l21}, the subgraph induced by $C$
contains at most $(d-1)c-{d \choose 2}$ blue edges and at most $(d-1)c-{d \choose 2}$ red edges.
Hence, this subgraph contributes to the sum of $b(v)$ at most $2(d-1)c-d(d-1)$ and to the sum of $r(v)$
at most $2(d-1)c-d(d-1)$.
Hence, the sum of $b(v)$ ($r(v)$) on the vertices of $B$ ($R$) must be at least
$b_c-2(d-1)c+d(d-1)$ ($r_c-2(d-1)c+d(d-1)$). It follows that:
$$
2(d-1)(b+c)-{d \choose 2} \geq \sum_{v \in B \cup C} b(v) \geq b_c+(b_c-2(d-1)c+d(d-1)),
$$
$$
2(d-1)(r+c)-{d \choose 2} \geq \sum_{v \in R \cup C} r(v) \geq r_c+(r_c-2(d-1)c+d(d-1)).
$$
Summing the two last inequalities we have:
$$
2(d-1)(b+r)-d(d-1)+4(d-1)c \geq (2k-4(d-1))c+2d(d-1).
$$
Thus, $r+b \geq (k-4d+4)c/(d-1)+3d/2$.
On the other hand $r+b+c \leq n$.
It follows that
$$
c \leq \frac{d-1}{k-3d+3}n-\frac{3d(d-1)}{2(k-3d+3)}, \qquad \frac{r+b}{2}+c \leq \frac{k-2d+2}{2(k-3d+3)}n-
\frac{3d(d-1)}{4(k-3d+3)}.
$$
It follows that there is either a red or a blue monochromatic $d$-subgraph of order at least
$$
\frac{k-4d+4}{2(k-3d+3)}n+\frac{3d(d-1)}{4(k-3d+3)}.
$$
\npf

\noindent
{\bf Proof of Theorem \ref{t1}, part (ii).}\,
It suffices to prove the theorem for $k=4d-4$.
We first create a specific graph $H$ on $n$ vertices.
Place the $n$ vertices in a sequence $\{v_1,\ldots,v_n\}$ and connect any two vertices whose distance is at most
$d-1$. Hence, all the vertices $\{v_d,\ldots,v_{n-d+1}\}$ have degree $2(d-1)$. The first $d$ and last $d$
vertices have smaller degree. To compensate for this we add the following ${d \choose 2}$ edges.
For all $i=1,\ldots,d-1$ and for all $j=i,\ldots,d-1$ we add the edge $(v_i,v_{jd+1})$.
Hence, if, say, $d=3$ we add $(v_1,v_4)$, $(v_1,v_7)$ and $(v_2,v_7)$.
Notice that these added edges are indeed new edges. The resulting graph $H$ has $n$ vertices and
$(k-1)n$ edges. Furthermore, all the vertices have degree $2(d-1)$ except for
$v_{jd+1}$ whose degree is $2(d-1)+j$ for $j=1,\ldots,d-1$ and $v_{n-d+1+j}$ whose degree is
$2(d-1)-j$ for $j=1,\ldots,d-1$. Also notice that any $d$-subgraph
of $H$ may only contain the vertices $\{v_1,\ldots,v_{d^2-d+1}\}$. Thus, the order of
any $d$-subgraph of $H$ is at most $d^2-d+1$. The crucial point to observe is that
the vertices of excess degree, namely $\{v_{d+1},v_{2d+1},\ldots,v_{d^2-d+1}\}$
form an independent set. Hence, for $n$ sufficiently large, $K_n$ contains two edge disjoint
copies of $H$ where in the second copy, the vertex playing the role of $v_{jd+1}$ plays the
role of the vertex $v_{n-d+1+j}$ in the first copy, for $j=1,\ldots,d-1$, and vice versa.
In other words, there exists a $4(d-1)$-regular graph with $n$ vertices, and a red-blue coloring
of it, such that the red subgraph and the blue subgraph are each isomorphic to $H$. In particular,
there is no monochromatic $d$-subgraph with more than $d^2-d+1$ vertices. \npf

\noindent
{\bf Proof of Theorem \ref{t2}.}\,
The theorem is trivial for $d=1$ so we assume $d \geq 2$.
It clearly suffices to prove the theorem for $n=(m+d)r$ where $m$ is an arbitrary
element of some fixed infinite {\em arithmetic} sequence whose difference and first element
are only functions of $d,k$ and $r$. Let $m$ be a positive integer such that
$$
y=m\frac{(d-1)(r-1)}{k-(r+1)(d-1)}
$$
is an integer.
Whenever necessary we shall assume $m$ is sufficiently large.
We shall create a graph with $n=(m+d)r$ vertices, minimum degree at least $k$, having an
$r$-coloring of its edges with no monochromatic subgraph larger than the value stated in the theorem.
Let $A_1,\ldots,A_r$ be pairwise disjoint sets of vertices of size $y$ each.
Let $B_1,\ldots,B_r$ be pairwise disjoint sets of vertices (also disjoint from the $A_i$)
of size $x=m+d-y$ each.
The vertex set of our graph is $\cup_{i=1}^r (A_i \cup B_i)$.
The edges of $G$ and their colors are defined as follows.
In each $B_i$ we place a graph of minimum degree at least $k-(r-1)(d-1)$,
and color its edges with the color $i$.
In each $A_i$ we place a $(d-1)$-degenerate graph with the maximum possible number of 
vertices of degree $2(d-1)$.  It is easy to show that such 
graphs exists with precisely $d$ vertices of degree $d-1$ and the rest are of degree $2(d-1)$.
Denote by $A_i'$ the $y-d$ vertices of $A_i$ with degree $2(d-1)$ in this subgraph and
put $A_i'' = A_i \setminus A_i'$.
Color its edges with the color $i$.
Now for each $j \neq i$ we place a bipartite graph whose sides are $A_i$ and $A_j \cup B_j$
and whose edges are colored $i$. The degree of all the vertices of $A_j \cup B_j$ in this subgraph
is $d-1$, the degrees of all the vertices of $A_i'$ are
at least $(k-(r+1)(d-1))/(r-1)$ and the degrees of all vertices of $A_i''$ in this subgraph are at least
$(k-r(d-1))/(r-1)$.
This can be done for $m$ sufficiently large since
$$
(y-d)\left\lceil \frac{k-(r+1)(d-1)}{r-1} \right\rceil+d \left\lceil \frac{k-r(d-1)}{r-1} \right\rceil \leq (d-1)(m+d).
$$
Notice that when $m$ is sufficiently large
we can place all of these $r(r-1)$ bipartite subgraphs such that their edge sets are pairwise disjoint
(an immediate consequence of Hall's Theorem).

By our construction, the minimum degree of the graph $G$ is at least $k$.
Furthermore, any monochromatic subgraph with minimum degree at least $d$ must be completely
placed within some $B_i$. It follows that
$$
f(n,k,d,r) \leq x = m+d-m\frac{(d-1)(r-1)}{k-(r+1)(d-1)}=n \frac{k-2r(d-1)}{r(k-(r+1)(d-1))}+C.
$$
\npf

\noindent
{\bf Proof of Theorem \ref{t3}.}\,
Suppose $n \geq R(4d+2k-5 ,4d+2k-5)$ where $R(a,b)$ is the usual Ramsey number.
Let $G$ be a a graph with $\delta(G)=n-k$ and fix a red-blue coloring of $G$.
Add edges to $G$ in order to obtain $K_n$. Note that at most $k-1$ new edges are incident with each vertex.
Color the new edges arbitrarily using the colors red and blue.
The obtained complete graph contains either a red or blue $K_{4d+2k-5}$ .
Deleting the new edges we get a monochromatic subgraph of $G$ on $4d+2k-5$ vertices and
minimum degree at least $4d+k-4 \geq 4d-3 \geq d$.
Now consider the largest monochromatic subgraph $Y$ with minimum degree at least $d$.
Hence, $|Y| \geq 4d+ 2k -5$. Assume, w.l.o.g., that $|Y|$ is red.
If $|Y| \leq n-2d-k+2$, then define $X$ to be a set of $2d+k-2$ vertices in $V \setminus Y$.
We call a vertex $y \in Y$ {\em bad} if it has $d$ ``red'' neighbors in $X$.
Let $B$ denote the subset of bad vertices in $Y$.
Since the number of red edges between $X$ and $B$ is at most $|X|(d-1)$ we have
$|B|d \leq |X|(d-1)$. Hence, $|B| < |X|=2d+k-2 \leq 4d+2k-5 \leq |Y|$. In particular, $|B| \leq 2d+k-3$.
Consider the bipartite blue graph on $X$ versus $Y \setminus B$. Its order is $|X|+|Y|-|B| > |Y|$.
Furthermore, we claim that it has minimum degree at least $d$. This is true because
each $y \in Y \setminus B$ has at least $|X|-(d-1)-(k-1) = d$ blue neighbors in $|X|$
and each vertex in $X$ is adjacent to at least $|Y| - |B| -(d-1)- (k-1) \geq 4d+2k-5 - (2d +k -3)-(d-1) -(k-1) = d$
vertices in $Y \setminus B$. Thus, $X \cup (Y \setminus B)$ contradicts the maximality of $Y$.
So, we must have $|Y| \geq  n-2d-k+3$, as required.
Clearly the value $n-2d-k+3$ is sharp for large $n$.
Take a red $K_{n-2d-k+3}$ on vertices $v_1,\ldots,v_{n-2d-k+3}$ and a blue $K_{2d+k-3}$
on vertices $u_1,\ldots,u_{2d+k-3}$. Put $A = \{v_1,\ldots, v_{2d+k-3} \}$. Connect with
$d-1$ blue edges the vertex $u_i$ to the vertices
$v_i ,\ldots,v_{i+d-2 (\bmod 2d+k-3)}$ , and connect with $d-1$ red edges the vertex
$u_i$ to the vertices $v_{i+d-1},\ldots, v_{i+2d -3 (\bmod 2d+k-3)}$.
There are no edges between $u_i$ and $v_{i+2d-2}, \ldots, v_{i+2d+k-4 (\bmod 2d+k-3)}$ .
The rest of the edges between the $u_i$ and $v_j$ for $j \geq 2d+k-2$ are colored blue.
It is easy to verify that this graph is $(n-k)$-regular and contain no blue nor red $d$-subgraph with
more than $n-2d-k+3$ vertices.
\npf

\section{Concluding remarks}
\begin{itemize}
\item
In the proof of Theorem \ref{t3} we assume $n \geq R(4d+2k-5 ,4d+2k-5)$ and hence $n$ is very large.
We can improve upon this to $n \geq \Theta(d+k)$ using the following argument.
Let $g(n,m,d,r)$ denote the largest integer $t$ such that in any $r$ coloring of a graph with $n$ vertices and $m$
edges there exists a monochromatic subgraph of order at least $t$ and minimum degree $d$.
\begin{prop}
\label{p31}
$$
g(n,m,d,r) \geq \sqrt{2\left(m - (d-1)n +{d \choose 2}\right)/r} \geq \sqrt {2m/r - 2dn/r}.
$$
\end{prop}
{\bf Proof.}\,
Suppose $G$ has $n$ vertices $m$ edges and the edges are $r$-colored.
Start deleting edge-disjoint monochromatic $d$-graphs as long as we can.
We begin with $m$ edges and when we stop we remain with at most $(d-1)n - {d \choose 2}$ edges.
Hence, there are at least $q=(m -(d-1)n +{d \choose 2})/r$ edges in one of the monochromatic $d$-graphs.
Thus, this monochromatic $d$-graph contains at least $\sqrt{2q}$ vertices as claimed.
Notice that this bound is rather tight for $d \leq \sqrt{2m/r}-1$.
Consider the $n$-vertex graph composed of $r$ vertex-disjoint copies of
$K_{\sqrt{2m/r}}$ and $n-\sqrt{2mr}$ isolated vertices (assume all numbers are integers, for simplicity).
Then, $e(G) \geq m$ and by coloring each of the $r$ large cliques with different colors
we get that any monochromatic $d$-subgraph has at most $\sqrt{2m/r}$ vertices.
\npf

Proposition \ref{p31} shown that in the proof of Theorem \ref{t3}
we can ensure an initial big monochromatic $d$-subgraph already when $n \geq 7(k+2d)/2= \Theta(d+k)$.

\item
In the case where $r \geq 3$ colors are considered and $k > 2r(d-1)$ is fixed, Theorem \ref{t2} supplies a linear
upper bound for $f(n,k,d,r)$. However, unlike the case where only two colors are used, we do not have a matching
lower bound. The following recursive argument supplies a linear lower bound in case $k=k(d)$ is
sufficiently large.
We may assume that $r$ is a power of $2$ as any lower bound for $r$ colors implies a lower bound for
less colors. Given an $r$-coloring of an $n$-vertex graph $G$, split the colors into two groups of $r/2$ colors each.
Now, using Theorem \ref{t1} we have a subgraph that uses only the colors of one of the groups,
and whose minimum degree is $x$, where $x$ is a parameter satisfying $k \geq 4x-3$.
The order of this subgraph is at least $n(k-4x+4)/(2(k-3x+3))$. Now we can use the recursion to
show that this $r/2$-colored linear subgraph has a linear order subgraph which is monochromatic.
$x$ is chosen so as to maximize the order of the final monochromatic subgraph.
For example, with $r=4$ we can take $x=4d-3$ and hence $k \geq 16d-15$.
For this choice of $x$ (which is optimal for this strategy) we get a monochromatic subgraph of order
at least
$$
n \frac{k-4(4d-3)+4)((4d-3) - 4d +4)}{(2(k- 3(4d-3)+3))( 2( (4d-3) - 3d +3))}=
n \frac{k - 16d  +16}{4d(k-12d+12)}.
$$

\item
Our theorems determine, up to a constant additive term, the value of $f(n,k,d,2)$ whenever $k$ or $n-k$ are fixed
and $n$ is sufficiently large. It may be interesting to establish precise values for all $k < n$.
Another possible path of research is the extension of the definition of $f(n,k,d,r)$ to $t$-uniform hypergraphs.
\end{itemize}

\end{document}